\documentclass[11pt]{amsart}      
\usepackage{amsmath,amsthm,amssymb}  
\usepackage{latexsym}        
\usepackage{euscript}  
\usepackage[frame,ps,dvips,matrix,arrow,curve,rotate]{xy}    

\newtheorem{thm}{Theorem}[section]
\newtheorem{defn}[thm]{Definition}
\newtheorem{prop}[thm]{Proposition}
\newtheorem{claim}[thm]{Claim}
\newtheorem{cor}[thm]{Corollary}
\newtheorem{lemma}[thm]{Lemma}

\theoremstyle{definition}  
\newtheorem{example}[thm]{Example}
\newtheorem{examples}[thm]{Examples}
\newtheorem{facts}[thm]{Facts}
\newtheorem{fact}[thm]{Fact}
\newtheorem{exercise}[thm]{Exercise}
\newtheorem{note}[thm]{Note}
\newtheorem{remark}[thm]{Remark}
\newtheorem{notation}[thm]{Notation}

\newcommand{\cat}[1]{{\EuScript #1}}
\newcommand{\cA}{\cat{A}}      
      
\newcommand{\cC}{\cat{C}}

\newcommand{\cE}{\cat{E}}
\newcommand{\cF}{\cat{F}}
\newcommand{\cG}{\cat{G}}
\newcommand{\cH}{\cat{H}}

\newcommand{\cK}{\cat{K}}

\newcommand{\cO}{\cat{O}}

\newcommand{\cU}{\cat{U}}
\newcommand{\cV}{\cat{V}}
\newcommand{\cW}{\cat{W}}

\newcommand{\bi}{\begin{itemize}}
\newcommand{\ei}{\end{itemize}}
\newcommand{\bc}{\begin{center}}
\newcommand{\ec}{\end{center}}
\newcommand{\ba}{\begin{array}}
\newcommand{\ea}{\end{array}}
\newcommand{\bt}{\begin{tabular}}
\newcommand{\et}{\end{tabular}}
\newcommand{\beqn}{\begin{equation}}
\newcommand{\eeqn}{\end{equation}}

\newcommand{\bp}{\begin{proof}}
\newcommand{\ep}{\end{proof}}
\newcommand{\blemma}{\begin{lemma}}
\newcommand{\elemma}{\end{lemma}}
\newcommand{\bprop}{\begin{prop}}
\newcommand{\eprop}{\end{prop}}
\newcommand{\bthm}{\begin{thm}}
\newcommand{\ethm}{\end{thm}}
\newcommand{\bremark}{\begin{remark}}
\newcommand{\eremark}{\end{remark}}
\newcommand{\bcor}{\begin{cor}}
\newcommand{\ecor}{\end{cor}}
\newcommand{\bnote}{\begin{note}}
\newcommand{\enote}{\end{note}}
\newcommand{\bdefn}{\begin{defn}}
\newcommand{\edefn}{\end{defn}}
\newcommand{\bexample}{\begin{example}}
\newcommand{\eexample}{\end{example}}
\newcommand{\bclaim}{\begin{claim}}
\newcommand{\eclaim}{\end{claim}}
\newcommand{\bexercise}{\begin{exercise}}
\newcommand{\eexercise}{\end{exercise}}
\newcommand{\bnotation}{\begin{notation}}
\newcommand{\enotation}{\end{notation}}
\newcommand{\bexamples}{\begin{examples}}
\newcommand{\eexamples}{\end{examples}}
\newcommand{\bfacts}{\begin{facts}}
\newcommand{\efacts}{\end{facts}}
\newcommand{\bfact}{\begin{fact}}
\newcommand{\efact}{\end{fact}}

\newcommand{\field}[1]  {\mathbb{#1}} 
\newcommand{\R}         {\field{R}}

\newcommand{\Z}         {\field{Z}}
\newcommand{\C}         {\field{C}}
\newcommand{\D}         {\field{D}}

\newcommand{\Q}         {\field{Q}}

\newcommand{\sP}        {\field{P}}           
\newcommand{\T}{\hat{T}}
\renewcommand{\H}{\hat{H}}

\renewcommand{\L}{\hat{L}}

\newcommand{\ho}{\textit{\textbf{H}}}
\newcommand{\ko}{\textit{\textbf{K}}}

\newcommand{\eps}{\epsilon}
\renewcommand{\aa}{\alpha}
\newcommand{\bb}{\beta}
\renewcommand{\gg}{\gamma}
\newcommand{\Gg}{\Gamma}
\newcommand{\dd}{\delta}

\renewcommand{\ll}{\lambda}
\newcommand{\Ll}{\Lambda}

\renewcommand{\ss}{\sigma}

\DeclareMathOperator{\Hom}{Hom}
\DeclareMathOperator{\im}{Im}  
\DeclareMathOperator{\Spec}{Spec}

\DeclareMathOperator{\dist}{dist}
\DeclareMathOperator{\Ell}{Ell}
\newcommand{\specm}{\Spec_m}

\newcommand{\ts}{\otimes}                       
\newcommand{\ra}{\rightarrow}                   
\newcommand{\lra}{\longrightarrow}              
\newcommand{\la}{\leftarrow}                    
\newcommand{\llra}[1]{\stackrel{#1}{\lra}}      

\newcommand{\iso}{\llra{\sim}}                  
\newcommand{\inc}{\hookrightarrow}              

\newcommand{\vs}{\vspace{2mm}}

\newcommand{\blank}{-}                          
\newcommand{\period}    {{\makebox[0pt][l]{\hspace{2pt} .}}}
\newcommand{\comma}     {{\makebox[0pt][l]{\hspace{2pt} ,}}}

\begin{document}

\title{Equivariant K-theory and equivariant cohomology}
  
\maketitle

\bc By \textsc{Ioanid Rosu}\\
with an appendix by \textsc{Allen Knutson} and \textsc{Ioanid Rosu}
\ec

\begin{abstract}
For $T$ an abelian compact Lie group, we give a description of
$T$-equivariant $K$-theory with complex coefficients in terms of
equivariant cohomology.  In the appendix we give applications
of this by extending results of Chang--Skjelbred and
Goresky--Kottwitz--MacPherson from equivariant cohomology to
equivariant $K$-theory.
\end{abstract}

\section{Introduction}

Let $T$ be an abelian compact Lie group, not necessarily connected.
Let $X$ be a compact $T$-equivariant manifold, or more generally a finite
$T$-$CW$ complex.  We denote by $H_T^*(X)$ the $T$-equivariant (Borel)
cohomology of $X$, as described in Atiyah and Bott~\cite{AB}, and by
$K_T^*(X)$ the $T$-equivariant $K$-theory of $X$, as described in
Segal~\cite{Seg}.  All the cohomology theories in this paper have
complex coefficients, unless otherwise noted.  For example,
$K_T^*(X)=K_T^*(X,\Z)\ts_\Z \C$.  Also, when $X$ is a point, we
write $K_T^*$ instead of $K_T^*(X)$, and similarly for $H_T^*$.

The goal of this paper is to describe $K_T^*(X)$ in terms of
$H_T^*(X)$.  When $T$ is the trivial group, this is easy:
it is a classical result that the Chern character
$ch:K^*(X) \ra H^*(X)$ is an isomorphism (in this case one only
needs to tensor with $\Q$).  In general however it is not true that
the equivariant version of the Chern character 
      $$ch_T:K_T^*(X) \ra H_T^{**}(X)$$
is an isomorphism.  (For the definition of $ch_T$ see
Lemma~\ref{cht-holo} and the discussion before it.)

The good news is that there is still a way in which one can describe
$K_T^*(X)$ in terms of $H_T^*(X)$.  Details will be given later, but
for now let us outline the main steps of this description.  Let
$\cC_T = \specm K_T^*$ be the complex algebraic group of the maximal
ideals of $K_T^*$.  The construction of $\cC_T$ is functorial in $T$,
and if $H\inc T$ is a compact subgroup of $T$, we can identify $\cC_H$
as a subgroup of $\cC_T$ via the map $\cC_H \ra \cC_T$.
If $\aa$ is a point of $\cC_T$, denote by $H(\aa)$ the smallest
compact subgroup $H$ of $T$ such that $\cC_H$ contains $\aa$.
Denote by $X^\aa=X^{H(\aa)}$, the subspace of $X$ fixed by all
elements of $H(\aa)$.  Denote by $\cO$ the sheaf of algebraic functions
on $\cC_T$, and by $\cO^h$ the sheaf of holomorphic functions.
Then we will define a sheaf, denoted by $\cK_T^*(X)$, whose stalk at a
point $\aa\in\cC_T$ is
     $$\cK_T^*(X)_\aa = \ho_T^*(X^\aa) \comma$$
where $\ho_T^*(\blank)$ is the extension of $H_T^*(\blank)$
by the ring of holomorphic germs at zero on $H_T^*$.
Moreover, the transition functions of $\cK_T^*(X)$ will be also
defined entirely using the equivariant cohomology of $X$.
Now, if $\Gg$ denotes the global sections functor, we will show that
there exists an isomorphism
  $$\Phi_T: K_T^*(X) \ts_{\Gg\cO} \Gg\cO^h \cong \Gg\cK_T^*(X) \period$$
This is the sense in which equivariant $K$-theory with complex
coefficients can be described in terms of equivariant cohomology. 

We should say a few words about the isomorphism $\Phi_T$.
If we denote by $ch_T$ the equivariant Chern character (to be
defined below), we will see that $\Phi_T$ is esentially a sheaf
version of $ch_T$.  However, $ch_T$ has to be twisted
(translated) in an appropriate sense, to take into account the
point $\aa$ over which the germ of $ch_T$ is taken.

One may call this a de Rham model for equivariant $K$-theory, but
it would be a slight misnomer, since we do not describe $K_T^*(X)$
at the level of cocycles.  As a matter of fact, we do more: we
describe the classes themselves, as sections in a sheaf build from
ordinary equivariant cohomology.  But, if one were really intent on
giving a de Rham model, one could use the sheaf model to define
a cocycle in $K$-theory as a collection of germs of ordinary
closed differential forms, and use a similar definition for
coboundaries.  We do not purse this avenue because it would only
obscure the purely topological nature of our description of
$K_T^*(X)$.

There were previous attempts to give a de Rham type of model for
$K_T^*(X)$.  The earliest version appeared in Baum, Brylinski and
MacPherson~\cite{BBM}.  The ideas were further developed in Block
and Getzler~\cite{BG}, and Duflo and Vergne~\cite{DV}.  In fact,
both the idea of describing equivariant $K$-theory as a sheaf and
twisting the equivariant Chern character are present in Duflo and
Vergne.  The problem in their paper is that they cannot prove the
Mayer--Vietoris property for their cohomology theory, because they
work with $\cC^\infty$ functions.  The advantage of our approach is
that we work with coherent analytic sheaves over the affine (Stein)
manifold $\cC_T$, and in this case the global section functor is
exact.

The present paper is inspired mainly by Grojnowski's preprint~\cite{Gr}.
In this seminal work, he uses ideas from the papers mentioned above
to define equivariant elliptic cohomology with complex coefficients.
His model starts with an elliptic curve $\cE$, and constructs for
every torus $T$ a complex variety $\cE_T$ and a coherent analytic
sheaf $\Ell_T^*(X)$ over $\cE_T$, whose stalk at each point is
defined in terms of equivariant cohomology.  The sheaf $\Ell_T^*(X)$
is then defined by Grojnowski to be the (``delocalized'') complex
$T$-equivariant cohomology of $X$.
Interestingly enough, equivariant $K$-theory is never explicitly
mentioned in Grojnowski's preprint, although he was most likely
aware that an analogous construction to that of $\Ell_T^*(X)$ should
lead to equivariant $K$-theory, if the elliptic curve $\cE$ is
replaced by the multiplicative group $\cC = \C \setminus \{0\}$.
The main contribution of the present paper is to do
exactly that: it starts with the multiplicative group $\cC$, out
of which it defines the base complex variety $\cC_T$, and constructs
a coherent analytic sheaf $\cK_T^*(X)$ over $\cC_T$.  Then, the ring
of global sections in $\cK_T^*(X)$ turns out to be a faithfully flat
extension of equivariant $K$-theory\footnote{Besides its
contribution to equivariant K-theory, one can regard the present
paper as giving a rigorous definition for Grojnowski's equivariant
elliptic cohomology: for this, it is enough to change the base
manifold $\cC_T$ to Grojnowsi's $\cE_T$.  See also Rosu~\cite{Ro}
for a definition of equivariant elliptic cohomology when $T=S^1$.}.  

A simple exercise shows that if one starts instead with the additive
group $\cA=\C$, the resulting sheaf is nothing else but ordinary
equivariant cohomology.  An important conclusion is that, when
working with complex coefficients, the difference between equivariant
cohomology, $K$-theory and elliptic cohomology stems mainly from the
fact that these theories are associated to different complex groups
of dimension one: the additive, the multiplicative, and the elliptic
groups, respectively.

The results of this paper can be extended in several directions.
First, we can describe $K_T^*(X)$ directly, instead of describing its
faithfully flat extension $K_T^*(X) \ts_{\Gg\cO} \Gg\cO^h$.  But in
order to do that, one needs to define algebraic sheaves $\cF_\aa$
(see Definition~\ref{local-sheaf}) instead of holomorphic ones.  And
this can only be done using completions, because the logarithm map is not
algebraic.  The construction therefore becomes more complicated, and
we decided to relegate it to another paper.  Second, if $G$ is a
nonabelian connected compact Lie group, and $T$ its maximal torus,
then $K_G^*(\blank)=K_T^*(\blank)^W$,
where $W$ is the Weyl group, so one can also describe $K_G^*(\blank)$
using Borel equivariant cohomology.  Third, we can prove a similar
result whenever the coefficient ring $R$ of the cohomology theories
involved is an algebra over $\Q$ adjoined the roots of unity.
We need $R$ to be a $\Q$-algebra because the logarithm map is only
defined over $\Q$, and we need to invert the roots of unity because we
want to split $R[z]/ \langle z^n-1\rangle$ into
a direct sum of $n$ copies of $R$.

While the details of the sheafifying process are somewhat technical,
in principle the construction allows one to infer some results in
equivariant $K$-theory from the corresponding ones in equivariant cohomology.
In the Appendix we give examples of this, extending results of
Chang--Skjelbred~\cite{CS} and Goresky--Kottwitz--MacPherson~\cite{GKM} 
from equivariant cohomology to equivariant $K$-theory.

\section{A sheaf--valued cohomology theory}

The purpose of this section is to define a sheaf valued $T$-equivariant
cohomology theory, which we denote by $\cK_T^*(\blank)$.  In the
next sections we are going to show that global sections of this sheaf
are essentially equivariant $K$-theory.  We already knew that
$K_T^*(X)$ can be regarded as a coherent sheaf over
$\cC_T = \specm K_T^*$ (because it is a $K_T^*$-module).  The novelty
is that $K_T^*(X)$ can be completely described in terms of \emph{ordinary}
equivariant cohomology (since we will show that $\cK_T^*(X)$ is).  Let us start with
a few definitions.\footnote{For a similar definition in the case of equivariant elliptic cohomology, see Rosu~\cite{Ro}.  The discussion there is only for $T=S^1$, but it generalizes easily with the same formalism as here.}.

\vs
\bc
{\sc 2.1. Definitions}
\ec
\vs

First we want a simpler description of $\cC_T=\specm K_T^*$.
Denote by $\T$
the Pontrjagin dual of $T$, i.e.\ $\T=\Hom(T,S^1)$.  For example,
if $T=(S^1)^p \times G$, where $G$ is a finite abelian group, then
$\T\cong\Z^p \times G$ (the isomorphism is not natural, however).
Denote by $\cC$ the multiplicative algebraic group $\C\setminus\{0\}$.
Although it might generate some confusion, we will use additive
notation for $\cC$ throughout the paper.  The following straighforward
lemma gives two alternate descriptions of $\cC_T$.

\bprop  \label{ct}
There is a natural isomorphism of algebraic varieties
   $$ \cC_T \cong \Hom_\Z (\T,\cC) \period$$
If $T$ is connected (i.e.\ a torus), and $\Ll_T$ is its integer lattice,
then there is a natural isomorphism
   $$ \cC_T \cong \Ll_T \ts_\Z \cC \period $$
\eprop

\bp
It is easy to see that $K_T^*=\C[\T]$, the group algebra of $\T$.
We define a map
   $$ \nu:\Hom_\Z (\T,\cC) \ra \specm\C[\T]=\cC_T $$
by noting that $\aa\in \Hom_\Z(\T,\cC)$ extends to a non-zero
$\C$-algebra map $\aa':\C[\T] \ra \C$.  We then take 
$\nu(\aa)=\ker(\aa')$, which is a maximal ideal of $\C[T]$.
Since the domain and codomain of $\nu$
both take products of groups to products of varieties, it
suffices to check that $\nu$ is an isomorphism when $T=S^1$
or $T=\Z_n$, which we leave to the reader.  

For the second statement, notice that when $T$ is a torus,
there is a natural isomorphism $\T \iso \Ll_T^*$.
\ep

Suppose $T\cong (S^1)^p\times G$, where $G$ is a finite abelian group.
Then one can apply the previous proposition to obtain that 
$\cC_T\cong \cC^p\times G$.  This formula shows that $\cC_T$ is
in fact an algebraic group (and it is the complexification of $T$).
Just as we did for $\cC$, in the rest of the paper we are going to use
{\it additive} notation for $\cC_T$.

\bdefn \label{exp}
Let $t_\C$ be the complexification of the Lie algebra of $T$.  Then
the exponential map $\exp:\C \ra \cC$ extends to a complex algebraic
group map
        $$\exp: t_\C \ra \cC_T \period$$
\edefn

\noindent
To be more precise, this map is the composite map
$t_\C = \Ll_T \ts_\Z \C \ra \Ll_T \ts_Z \cC = \cC_{T^0} \ra \cC_{T}$,
where $T^0$ is the connected component of $T$ containing the identity.
(For the identification $\Ll_T \ts_Z \cC = \cC_{T^0}$, use the
previous proposition.)  Note that when $T$ is connected, the
exponential map is surjective.

\bdefn \label{small}
We call a neighborhood $U$ of zero in $\cC_T$ ``small'' if the above
defined exponential map, $\exp:t_\C \ra \cC_T$, has a local inverse on $U$.
We call a neigborhood $V$ of zero in $t_\C$ ``small'' if it is of the
form $\exp^{-1}(U)$, for $U$ a small neighborhood of zero in $\cC_T$.
\edefn

Let $\cA$ be a collection of compact subgroups of $T$.  Define a relation
on $\cC_T$ as follows: $\aa \prec_\cA \bb$ if, for any $H\in\cA$,
$\bb\in\cC_H$ implies $\aa\in\cC_H$.  The relation $\prec_\cA$ is reflexive
and transitive, but not antisymmetric, so it is not an order relation.
When it is clear what $\cA$ is, we will omit it and write simply
$\aa\prec\bb$.  The next definition singles out a special class of
open covers of $\cC_T$, called adapted covers.  These will be used
below in the definition of the sheaf $\cK_T^*(X)$.

\bdefn \label{adapted-defn}
Let $\cA$ be a collection of compact subgroups of $T$, and let
$\cU=(U_\aa)_{\aa\in\cC_T}$ be an open cover indexed by the points of
$\cC_T$.  Then $\cU$ is called ``adapted to $\cA$'' if it satisfies
the following conditions:
\bi
\item[1.] $\aa\in U_\aa$, and $U_\aa-\aa$ is small.
\item[2.] If $U_\aa \cap U_\bb \neq \emptyset$, then either $\aa\prec\bb$
  or $\bb\prec\aa$.
\item[3.] If $\aa\prec\bb$, and for some $H\in\cA$ $\aa\in\cC_H$ but
$\bb\notin\cC_H$, then $U_\bb \cap \cC_H = \emptyset$.
\item[4.] If $U_\aa \cap U_\bb \neq \emptyset$, and both $\aa$ and $\bb$
  belong to $\cC_H$ for some $H\in\cA$, then $\aa$ and $\bb$ belong to the
  same connected component of $\cC_H$.
\ei
\edefn

\bprop \label{adapted}
If $\cA$ is a finite collection of compact subgroups of $T$, then there
exists a cover $\cU$ of $\cC_T$ adapted to $\cA$.  Any refinement of
$\cU$ is still adapted.
\eprop

\bp
Define $\cH=\{\cC_H \; | \; H\in\cA \}$, and $\cH^0=$ the set
of all connected components of the elements in $\cH$.  Put a metric on
$\cC_T$ which yields its usual topology.  Denote this metric by
``$\dist$''.

Let $\aa\in\cC_T$.  If $\aa\in C$ for all $C\in\cH^0$ (this is possible
only when $T$ is connected), then choose $U_\aa$ such that $\aa\in U_\aa$,
and $U_\aa-\aa$ is small.  If, on the contrary, there exists a connected
component $C\in\cH^0$ such that $\aa\notin C$, then take $U_\aa$ a ball of
center $\aa$ and radius $d$, with
   $$d<\tfrac12 \min_{D\in\cH^0,\aa\notin D} \dist(\aa,D) \comma$$
and such that $U_\aa-\aa$ is small.

We show that $\cU=(U_\aa)_{\aa\in\cC_T}$ is adapted:  Condition 1
is trivially satifsfied.  To prove Condition 2,
let $\aa$ and $\bb$ be such that $U_\aa \cap U_\bb \neq \emptyset$.
Suppose we have neither $\aa\prec\bb$, nor $\bb\prec\aa$.  Then by the
definition of $\prec\,$ there exist two compact subgroups $K$ and $L$ of $T$ such
that $\aa\in\cC_K\setminus\cC_L$ and $\bb\in\cC_L\setminus\cC_K$.
But from the definition of $U_\aa$ it follows that $U_\aa$ is a ball 
of center $\aa$ and radius 
$d<\tfrac12\dist(\aa,\cC_L)\leq \tfrac12\dist(\aa,\bb)$.  Similarly, $U_\bb$
is a ball of center $\bb$ and radius less than $\tfrac12\dist(\aa,\bb)$,
so $U_\aa$ and $U_\bb$ cannot possibly intersect, contradiction.

Condition 3 is obviously satisfied, by construction.

Finally, to show Condition 4, let $\aa,\bb \in \cC_H$ be such that
$U_\aa \cap U_\bb \neq \emptyset$.  Suppose $\aa$ and $\bb$ belong to
different connected components of $\cC_H$.  Then by the same type of
reasoning as above, it follows that the radii of $U_\aa$ and $U_\bb$
are smaller than $\tfrac12\dist(\aa,\bb)$, so $U_\aa$ and $U_\bb$ cannot
possibly intersect, which again leads to a contradiction.
\ep

Let $\aa\in\cC_T$.  The construction of $\cC_T$ is functorial, so if
$H$ is any compact subgroup of $T$, we get an inclusion map $\cC_H \ra \cC_T$.
If $\aa \in \im(\cC_H \ra \cC_T)$, we say that $\aa \in \cC_H$.
For $\aa\in\cC_H$, denote by $H(\aa)$ the smallest compact subgroup $H$ of
$T$ such that $\aa\in\cC_H$.  The fact that there exists a smallest
$H$ such that $\aa\in\cC_H$ is implied by the formula
$\cC_K \cap \cC_L = \cC_{K\cap L}$, which follows from an easy diagram
chase.  This also implies that $H(\aa)$ is the intersection of all
compact subgroups $H$ such that $\aa\in\cC_H$.  Another immediate
consequence is the following formula, which will be useful later:
\begin{equation} \label{haa}
  H(\aa) \subseteq K \iff \aa\in\cC_K \period
\end{equation}

Let $X$ be a space with a $T$-action.  If $K$ is a compact subgroup
of $T$, denote by $X^K\subseteq X$ the subspace of points fixed
by $K$.  Also, if $\aa\in\cC_T$, define 
           $$X^\aa = X^{H(\aa)} \period$$
Now we want to define the notion of a cover adapted to a finite
$T$-$CW$ complex.  So let $X$ be a finite $T$-$CW$ complex.  We know
that there exists a finite collection $\cA=(H_i)_i$ of compact subgroups
of $T$ such that $X$ has an equivariant cell decomposition of the form
$X = \bigcup_i \D^{n_i}\times (T/H_i)$.  Here we denoted by
$\D^n$ the open disk in dimension $n$, and by $\D^0$ a point.
The group $T$ acts trivially on $\D^{n_i}$, and by left multiplication
on $T/H_i$.
Notice that if $K$ is a compact subgroup of $T$, then the subcomplex of
$X$ fixed by $K$ has a decomposition of the form
$X^K=\bigcup_{i:K\subseteq H_i} \D^{n_i}\times (T/H_i)$.
Taking $K=H(\aa)$, we get
\begin{equation} \label{xaa}
   X^\aa= \bigcup_{i:\aa\in\cC_{H_i}} \D^{n_i}\times (T/H_i)\period
\end{equation}
We say that the cover $\cU=(U_\aa)_{\aa\in\cC_T}$ of $\cC_T$ is
``adapted to $X$'' if $\cU$ is adapted to
the collection $\cA$ of the isotropy groups $H_i$ appearing in
a $T$-equivariant cell decomposition of $X$.  Since this collection
is finite, Proposition~\ref{adapted} implies that there always
exists a cover adapted to $X$.

Next we discuss a few useful results in equivariant cohomology.
We start with a well-known proposition which says that the ring
of coefficients of (complex) $T$-equivariant cohomology is the
polynomial algebra on $t_\C$, the complex Lie algebra of $T$.

\bprop \label{polynomial}
If $T$ is an abelian compact Lie group, there is a natural isomorphism
  $$ S(t_\C^*) \iso H_T^* \comma$$
where $S(\blank)$ denotes the symmetric algebra, and $t_\C^*$
is the dual of $t_\C$.
\eprop

\bp
$\T=\Hom(T,S^1)$ is the group of irreducible characters of
$T$, so for $\ll\in\T$ consider the complex vector bundle
       $$V_\ll= ET \times_T \C \comma$$
over the classifying space $BT$, where the map $T\ra\C$ is given
by $\ll$.  Then the first Chern class of $V_\ll$ gives a
natural isomorphism $c_1: \T \ra H^2(BT,\Z)$.

If $T$ is a torus, we saw in Proposition~\ref{ct} that $\T$
can be identified with the dual of the integer lattice $\Ll_T^*$,
so by tensoring the map $c_1$ with $\C$, we get the natural
isomorphism $c_1:t_\C^*=\Ll_T^*\ts\C \ra H^2(BT,\C)$.  Taking
symmetric products, we get the desired isomorphism.

If $T$ is a general (non-connected) compact abelian Lie group,
the isomorphism still holds, since both domain and codomain
depend only on the connected component of $T$ containing the
identity, which is a torus. 
\ep

We now define an algebra homomorphism 
   $$h: H_T^* \ra \cO^h_{t_\C,0}$$
by taking a polynomial in $H_T^*=S(t_\C^*)$ and sending it to
its germ at zero.  The map is injective, so we can consider
$H_T^*$ as a subring of $\cO^h_{t_\C,0}$.
Let $V$ be a small neighborhood of zero in $t_\C$.  Then, since
the ring $H_T^*\subset \cO^h_{t_\C,0}$ consists of the germs of
global holomorphic functions on $t_\C$, the map $h$ factors through
the inclusion $\cO^h_{t_\C}(V) \inc \cO^h_{t_\C,0}$, so we can 
define a map, also denoted by $h$,
   $$ h: H_T^* \ra \cO^h_{t_\C}(V) \period $$
Let $U$ be a small neighborhood of zero in $\cC_T$, and let
$V=\exp^{-1}(U)$, where $\exp:t_\C \ra \cC_T$ is the exponential
map.  Via the exponential, we have the following identifications:
$\cO^h_{\cC_T,0}\iso \cO^h_{t_\C,0}$ and
$\cO^h_{\cC_T}(U)\iso \cO^h_{\cC_T}(V)$.  

\bdefn \label{h}
Let $U$ be a small neighborhood of zero in $\cC_T$.
Via the identifications above, we define the following two maps,
and denote them also by $h$ (the second one is the correstriction
of the first):
  $$ h:H_T^* \ra \cO^h_{\cC_T,0} \quad\mbox { and } \quad
      h:H_T^* \ra \cO^h_{\cC_T}(U) \period $$
\edefn

Now we define a few cohomology theories that we are going
to use throughout the paper.  Let $X$ be a finite $T$-$CW$ complex.
We define the holomorphic $T$-equivariant cohomology of $X$ to be
  $$ \ho_T^*(X) =  H_T^*(X) \ts_{H_T^*} \cO^h_{\cC_T,0} \comma$$
where the map $h:H_T^* \ra \cO^h_{\cC_T,0}$ is given in
Definition~\ref{h}.  It is indeed a cohomology theory, because
by Proposition~\ref{flat-completion} the map
$H_T^* \ra \cO^h_{\cC_T,0}$ is flat.  The theory is not $\Z$-graded
anymore; however, it can be thought of as $\Z/2$-graded, by its
even and odd part.

Let $H_T^{**}(X)$ be the completion of $H_T^*(X)$ with respect
to the augmentation ideal $I=\ker(H_T^*\ra\C)$.  Since $H_T^*(X)$
is a finitely generated module over the Noetherian ring $H_T^*$,
a simple result on completions (see for example Matsumura~\cite{Ma},
Theorem 55) implies that
   $$H_T^{**}(X) \cong H_T^*(X) \ts_{H_T^*} H_T^{**} \period$$
The ring $H_T^*$ is a polynomial ring, so we have the following
well-known results from algebra (they are sometimes called {\it GAGA}
results, since they first appeared in Serre's {\it GAGA}~\cite{Ser}).

\bprop \label{flat-completion}
$\ho_T^*=\cO^h_{\cC_T,0}$ and $H_T^{**}$ are flat
over $H_T^*$.  If $U$ is a small neighborhood of zero, then
$\cO^h_{\cC_T}(U)$ is flat over $H_T^*$.
\eprop

\bp
All the results we mention in this proof are from Matsumura~\cite{Ma}.
Identify the complex Lie algebra $t_\C$ with $X=\C^n$.  We denote
by $\cO$ the algebraic structure sheaf of $\C^n$, and by $\cO^h$ the
analytic structure sheaf.  Let $\cO^\wedge_0$ be the completion of
the local ring $\cO_0$ with respect to its maximal ideal.  It is sufficient
to show that the natural maps $\cO(X) \ra \cO^h_0$ and
$\cO(X) \ra \cO^\wedge_0$ are flat, and that for any $U$ open in $X$, the
map $\cO(X) \ra \cO^h(U)$ is flat.  We start by noticing that $\cO^\wedge_0$
is the completion of $\cO(X)$ with respect to its maximal ideal at zero,
so by Corollary 1 of Theorem 55, we know that $\cO^\wedge_0$ is flat over
$\cO(X)$.  The completion of $\cO^h_0$ with respect to its maximal ideal
is also $\cO^\wedge_0$, so $\cO^h_0 \ra \cO^\wedge_0$ is flat.  It is in
fact faithfully flat, because it is local: see Theorem 3 (4.D).
Now one can check directly by the definition of flatness that having
$\cO(X) \ra \cO^\wedge_0$ flat and $\cO^h_0 \ra \cO^\wedge_0$ faithfully
flat implies that $\cO(X) \ra \cO^h_0$ is flat.  Notice also that the
same proof can be used to show that $\cO_0 \ra \cO^h_0$ is flat.

Next let $U\subseteq X$ be an open set.  By the local characterization
of flatness, Theorem (3.J), in order to show that $\cO(X) \ra \cO^h(U)$
is flat, we have to show that for any $x\in U$ the natural map
$\cO_x \ra \cO^h_x$ is flat.  But we have already shown this when
$x=0$, and the proof for general $x$ is the same.

Now in order to prove the proposition, just transfer the results
we have proved via the exponential map, $\exp:\C^n=t_\C \ra \cC_T$.
This is where we need $U$ small.
\ep

In particular, it follows that $H_T^*(X)$ and $\ho_T^*(X)$ can be
regarded as subrings of $H_T^{**}(X)$.

\vs
\bc
{\sc 2.2. Construction of} $\cK_T^*$
\ec
\vs

Let $X$ be a finite $T$-$CW$ complex.  Fix $\cU$ a cover
adapted to $X$, which exists because of Proposition~\ref{adapted}.
We are going to define a sheaf $\cF=\cF_\cU$
over $\cC_T$ whose stalk at $\aa\in\cC_T$ is isomorphic
to $\ho_T^*(X^\aa)$.  Recall that in order to give a sheaf
$\cF$ over a topological space, it is enough to give an open cover
$(U_\aa)_\aa$ of that space, and a sheaf $\cF_\aa$ on each $U_\aa$
together with isomorphisms of sheaves 
$\phi_{\aa\bb}: \cF_{\aa\mid_{U_\aa \cap U_\bb}} \lra 
                        \cF_{\bb\mid_{U_\aa \cap U_\bb}}$,
such that $\phi_{\aa\aa}$ is the identity function, and the cocycle
condition $\phi_{\bb\gg}\phi_{\aa\bb}=\phi_{\aa\gg}$ is satisfied on
$U_\aa \cap U_\bb \cap U_\gg$.  If $U$ is a subset of
$\cC_T$, denote by $U+\aa=\{x+\aa \; | \; x\in U\}$ the translation
of $U$ by $\aa$.

\bdefn \label{local-sheaf}
Define a presheaf $\cF_\aa$ on $U_\aa$ by taking, for any open
$U\subseteq U_\aa$,
 $$\cF_\aa(U) = H_T^*(X^\aa) \ts_{H_T^*} \cO^h_{\cC_T}(U-\aa) \comma$$
where the restriction maps are induced from those of $\cO^h_{\cC_T}$.
The ring map $h: H_T^* \ra \cO^h_{\cC_T}(U-\aa)$ is given in
Definition~\ref{h}.
\edefn

\bprop
$\cF_\aa$ is a coherent sheaf of $\cO^h_{\cC_T}$-modules.
\eprop

\bp
First we show that $\cF_\aa$ is a sheaf of $H_T^*$-modules.  If
$(U_i)_i$ is an open cover of a topological space $Y$, denote by
$U_{ij}=U_i\cap U_j$, etc.  Then a presheaf $\cG$ is a sheaf if
and only if for any $m>0$ and any finite cover $(U_i)_{i=1\ldots m}$
the following sequence is exact 
  $$ 0 \lra \cG(Y) \llra{r_0} \prod_i \cG(U_i) \llra{r_1}
   \prod_{i<j} \cG(U_{ij}) \llra{r_2} \cdots \lra
     \cG(U_{1\ldots m}) \lra 0 \comma$$
where $r_0(s)_i=s_{\,\mid U_i}$, 
$r_1\bigl(\prod_i s_i\bigr)_{i,j} = 
  s_{i\;\mid U_{ij}}-s_{j\;\mid U_{ij}}$,
etc.  Denote this sequence by $M(\cG)^\bullet$.  Since 
$\cO^h_{\cC_T}$ is a sheaf of $H_T^*$-modules on $Y=U_\aa$, the
sequence $M(\cO^h_{\cC_T})^\bullet$ is exact.  Now, by
Proposition~\ref{flat-completion}, each $H_T^*$-module in
the sequence $M(\cO^h_{\cC_T})^\bullet$ is flat.  Therefore,
$M(\cO^h_{\cC_T})^\bullet$ can be thought as a flat 
resolution of the last term to the right, $\cO^h_{\cC_T}(U_{1\ldots m})$.
Using the Flat Resolution Lemma (see Weibel~\cite{We}, Lemma 3.2.8),
one can easily show that by tensoring $M(\cO^h_{\cC_T})^\bullet$
with any $H_T^*$-module, we also get an exact sequence.  In
particular, by tensoring with $H_T^*(X^\aa)$ we obtain
$M(\cF_\aa)^\bullet$, which therefore is exact.  But that
means that $\cF_\aa$ is a sheaf (of $H_T^*$-modules).

To show that $\cF_\aa$ is a sheaf of $\cO^h_{\cC_T}$-modules,
if $U \subseteq U_\aa$, we need  an action of
$f\in\cO^h_{\cC_T}(U)$ on $\cF_\aa(U)$.  The translation map
$t_\aa: U-\aa \ra U$, which takes $x$ to $x+\aa$ gives a map
$t_\aa^*:\cO^h_{\cC_T}(U) \ra \cO^h_{\cC_T}(U-\aa)$, which sends
$f(u)$ to $f(u+\aa)$.  Then we take the result of the action
of $f\in\cO^h_{\cC_T}(U)$ on 
$\mu\ts g \in \cF_\aa(U) = 
        H_T^*(X^\aa) \ts_{H_T^*} \cO^h_{\cC_T}(U-\aa)$
to be $\mu\ts (t_\aa^* f\cdot g)$.  It is an easy exercise to show
that $\cF_\aa$ is coherent: one uses that $H_T^*(X^\aa)$ is a
finitely generated $H_T^*$-module.
\ep

Next, for the second part of the definition of $\cF$, 
we glue the different sheaves $\cF_\aa$.
Let $\aa,\bb\in\cC_T$ such that $U_\aa\cap U_\bb\neq\emptyset$.
Since the cover $\cU$ is adapted to $\cA$, Condition 2 from
Definition~\ref{adapted-defn} implies
that either $\aa\prec\bb$ or $\bb\prec\aa$.  Without loss of generality
we can assume that $\aa\prec\bb$.  Then it is clear that
$H(\aa)\subseteq H(\bb)$, and hence $X^\bb \subseteq X^\aa$.
Denote by $i:X^\bb \ra X^\aa$ the inclusion map.  Now we want to
investigate conditions in which the restriction map
$i^*: H_T^*(X^\aa) \ra H_T^*(X^\bb)$ becomes an isomorphism.
The next proposition is esentially the Localization Theorem
of Borel and Hsiang (see Atiyah and Bott~\cite{AB}).

\bprop \label{phi-localization}
Let $X$ be a finite $T$-$CW$ complex, and $\cU=(U_\aa)_{\aa\in\cC_T}$
a cover adapted to $X$.  Let $\aa\prec\bb$ be two points in $\cC_T$.
Let $U\subseteq U_\aa \cap U_\bb$, and assume $U_\aa \cap U_\bb$ nonempty.
Then the restriction map
  $$i^* \ts 1 : H_T^*(X^\aa) \ts_{H_T^*} \cO^h_{\cC_T}(U-\aa) \ra
         H_T^*(X^\bb) \ts_{H_T^*} \cO^h_{\cC_T}(U-\aa)$$
is an isomorphism.
\eprop

\bp
The case when $X^\aa=X^\bb$ is trivial, so we assume the contrary.
Then there exists a compact subgroup $L$ of $T$ such that $\aa\in\cC_L$
and $\bb\notin\cC_L$.  Then Condition 3 from 
Definition~\ref{adapted-defn}
implies that $U_\bb \cap \cC_L = \emptyset$.  But
$\cC_L = \aa+\cC_L$, so we deduce that $(U-\aa) \cap \cC_L = \emptyset$.
This means that for every $x\in U-\aa$ we have $x\notin\cC_L$.

Now we will prove that whenever $x\notin\cC_L$, $i^*$ becomes an
isomorphism if tensored with $\cO^h_{\cC_T,x}$.  This is enough
for us, because we can regard the statement to be proved as a
statement about sheaves, and we are proving the isomorphism
stalkwise.

Now notice that Equation~(\ref{xaa}) implies that $X^\aa$ is built
from $X^\bb$ by adding cells of the form $\D^n\times (T/L)$, with
$\aa\in\cC_L$ and $\bb\notin\cC_L$.  So it is enough to show that,
when tensored with $\cO^h_{\cC_T,x}$ over $H_T^*$, the ring
$H_T^*\bigl(\D^n\times (T/L),S^{n-1}\times (T/L)\bigr)$ becomes zero.
Or, equivalently, it is enough to have $H_T^*(T/L)=H_L^*$
become zero after tensoring with $\cO^h_{\cC_T,x}$.  We therefore
proceed to show that $H_L^* \ts_{H_T^*} \cO^h_{\cC_T,x}$ is the zero
algebra.

Via the correspondence $\log:t_\C \ra \cC_T$, it is enough
to show that if $l_\C$ is the complexification of the Lie algebra
of $L$, and $x\notin\ l_\C$, then $H_L^* \ts_{H_T^*} \cO^h_{t_\C,x}=0$.
Denote by $I_L$ the kernel of the surjective map of algebras
$H_T^*\ra H_L^*$.  It is clear that $I_L$ kills $H_L^*$ as an
$H_T^*$-module.  The vector space $l_\C$ is the vanishing set of $I_L$,
so since $x\notin l_\C$, there exists a polynomial $p$ in $I_L$ such that
$p(x) \neq 0$.  This means that $p$ is invertible in $\cO^h_{t_\C,x}$,
hence $1\ts p$ is invertible in $H_L^* \ts_{H_T^*} \cO^h_{t_\C,x}$.
By the balancing property of the tensor product, we know that $1\ts p=0$,
because $p$ goes to zero via the map $H_T^*\ra H_L^*$.  So we found an
element in $H_L^* \ts_{H_T^*} \cO^h_{t_\C,x}$ which is invertible and zero
at the same time.  This can only happen if $H_L^* \ts_{H_T^*} \cO^h_{t_\C,x}$
is the zero algebra.
\ep

Let us now consider the translation map $t_{\bb-\aa}: U-\bb \ra U-\aa$.
This induces a map 
$t_{\bb-\aa}^* : \cO^h_{\cC_T}(U-\aa) \ra \cO^h_{\cC_T}(U-\bb)$.  Eventually
we want to produce a translation map
$H_T^*(X^\bb) \ts_{H_T^*} \cO^h_{\cC_T}(U-\aa) \ra
    H_T^*(X^\bb) \ts_{H_T^*} \cO^h_{\cC_T}(U-\bb)$.  But the problem is
that $t_{\bb-\aa}^*$ is not a map of $H_T^*$-modules, so we cannot
simply tensor it with the identity map of $H_T^*(X^\bb)$.  However,
we hope to find some other group $K$ so that $t_{\bb-\aa}^*$ is a map
of $H_K^*$-modules.

\blemma \label{translation-k}
Let $X$ be a finite $T$-$CW$ complex, and $\cU=(U_\aa)_{\aa\in\cC_T}$
a cover adapted to $X$.  Let $\aa\prec\bb$ be two points in $\cC_T$ and
$U\subseteq U_\aa \cap U_\bb$, and assume $U_\aa \cap U_\bb$ nonempty.
Let $H=H(\bb)$ and $K=T/H$.  Then the translation map
  $$t_{\bb-\aa}^* : \cO^h_{\cC_T}(U-\aa) \ra \cO^h_{\cC_T}(U-\bb)$$
is a map of $H_K^*$-algebras.  Here the $H_T^*$-algebras
$\cO^h_{\cC_T}(U-\aa)$ and $\cO^h_{\cC_T}(U-\bb)$ are regarded
as $H_K^*$-algebras via the natural algebra map $H_K^* \ra H_T^*$
induced by the quotient map $T\ra K$.
\elemma

\bp
Since $\bb\in\cC_H$, it follows that $\aa\in\cC_H$.
Let $H^0$ be the connected component of $H$ containing the identity.
Because both $\aa$ and $\bb$ belong to $\cC_H$, Condition 4 from
Definition~\ref{adapted-defn} implies that $\aa$ and $\bb$ belong
to the same connected component of $\cC_H$, or equivalently that
$\bb-\aa\in \cC_{H^0}$.  Let $h_\C$ be the complexification of the
Lie algebra of $H$.  The exponential map $\exp:h_\C \ra \cC_{H^0}$
is a surjective group homomorphism, because $H^0$ is a torus.
Consider an element in the coset $\exp^{-1}(\bb-\aa) \subset h_\C$.
Denote it also by $\bb-\aa$.  Let $\Gg$ be the global sections
functor.  Using the exponential map, it is enough to show that
   $$t_{\bb-\aa}^* : \Gg\cO^h_{t_\C} \ra \Gg\cO^h_{t_\C}$$
is a map of $H_K^*$-modules.  This is equivalent to showing that
$t_{\bb-\aa}^*(p)=p$ for all $p\in H_K^*$.  But this is easy,
since we can identify $H_K^*$ with the ring of polynomial functions
$f(u)$ in $H_T^*$ such that $f(u+a)=f(u)$ for any $a$ in the
complexification of the Lie algebra of $K$.
\ep

Here is a well-known result in equivariant cohomology theory
(see for example Atiyah and Bott~\cite{AB} for ordinary cohomology
and Segal~\cite{Seg} for $K$-theory).

\bprop \label{fixed}
Let $h_T^*(\blank)$ be either ordinary equivariant cohomology,
or equivariant $K$-theory.
Let $Z$ be a finite $T$-$CW$ complex which is fixed by a compact subgroup
$H$ of $T$.  Let $K=T/H$.  Then there exists a natural isomorphism
   $$ h_K^*(Z) \ts_{h_K^*} h_T^* \iso h_T^*(Z)\period$$
\eprop

The last step in the construction of the gluing maps $\phi_{\aa\bb}$
is the following proposition:

\bprop \label{phi-translation}
Under the hypotheses of Lemma~\ref{translation-k}, the translation
map $t_{\bb-\aa}^*$ extends to a natural isomorphism
  $$\tau_{\bb-\aa}^* : H_T^*(X^\bb) \ts_{H_T^*} \cO^h_{\cC_T}(U-\aa) \ra
    H_T^*(X^\bb) \ts_{H_T^*} \cO^h_{\cC_T}(U-\bb) \period $$
\eprop

\bp
$X^\bb$ is fixed by $H=H(\bb)$, so by Proposition~\ref{fixed} we know
that $H^*_T(X^\bb) \ts_{H_T^*} \cO^h_{\cC_T}(U-\aa) \cong
   H^*_K(X^\bb) \ts_{H_K^*} \cO^h_{\cC_T}(U-\aa)$, with $K=T/H$.
Therefore it is enough to define a map
   $$ \phi:H_K^*(X^\bb) \ts_{H_K^*} \cO^h_{\cC_T}(U-\aa) \ra
    H_K^*(X^\bb) \ts_{H_K^*} \cO^h_{\cC_T}(U-\bb) \period $$
But  $t_{\bb-\aa}^*$ was shown in
Lemma~\ref{translation-k} to be a map of $H_K^*$-algebras,
so we can just take $\phi=1\ts t_{\bb-\aa}^*$.

More formally, we define $\tau_{\bb-\aa}^*$ as the composition
of the following isomorphisms:
 $H^*_T(X^\bb) \ts_{H_T^*} \cO^h_{\cC_T}(U-\aa) \ra
  H^*_K(X^\bb) \ts_{H_K^*} \cO^h_{\cC_T}(U-\aa) \llra{1\ts t_{\bb-\aa}^*}
  H^*_K(X^\bb) \ts_{H_K^*} \cO^h_{\cC_T}(U-\bb) \ra
  H_T^*(X^\bb) \ts_{H_T^*} \cO^h_{\cC_T}(U-\bb)$.
\ep

\bdefn \label{gluing-map}
Under the same hypotheses as Lemma~\ref{translation-k}, define
$\phi_{\aa\bb}:\cF_\aa(U) \ra \cF_\bb(U)$ as the composite map
$H_T^*(X^\aa) \ts_{H_T^*} \cO^h_{\cC_T}(U-\aa) \llra{i^*\ts 1}
   H_T^*(X^\bb) \ts_{H_T^*} \cO^h_{\cC_T}(U-\aa) \llra{\tau_{\bb-\aa}^*}
   H_T^*(X^\bb) \ts_{H_T^*} \cO^h_{\cC_T}(U-\bb)$,
where the second map was defined above.
\edefn

\bprop
$\phi_{\aa\bb}$ is an isomorphism.
\eprop

\bp
The first map is an isomorphism by Proposition~\ref{phi-localization},
and the second map is an isomorphism by Proposition~\ref{phi-translation}.
\ep

\bprop \label{phi-cocycle}
The maps $\phi_{\aa\bb}$ verify the cocycle condition, i.e.\
$\phi_{\bb\gg}\circ\phi_{\aa\bb}=\phi_{\aa\gg}$, on every
triple intersection $U_\aa\cap U_\bb\cap U_\gg\neq\emptyset$.
\eprop

\bp
Without loss of generality, we can assume $\aa\prec\bb\prec\gg$.
We need to show that the following diagram is commutative:
$$\xymatrix{
  H_T^*(X^\aa) \ts_{H_T^*} \cO^h_{\cC_T}(U-\aa) \ar[d]^{i^*\ts 1} & \\
  H_T^*(X^\bb) \ts_{H_T^*} \cO^h_{\cC_T}(U-\aa)
              \ar[d]^{\tau_{\bb-\aa}^*} \ar[dr]^{i^*\ts 1}  &    \\
  H_T^*(X^\bb) \ts_{H_T^*} \cO^h_{\cC_T}(U-\bb)  \ar[d]^{i^*\ts 1}  &
    H_T^*(X^\gg) \ts_{H_T^*} \cO^h_{\cC_T}(U-\aa) \ar[dl]^{\tau_{\bb-\aa}^*}\\
  H_T^*(X^\gg) \ts_{H_T^*} \cO^h_{\cC_T}(U-\bb) \ar[d]^{\tau_{\gg-\bb}^*} & \\
  H_T^*(X^\gg) \ts_{H_T^*} \cO^h_{\cC_T}(U-\gg) &
             }$$
In this diagram the only map which has not been defined is
$\tau_{\bb-\aa}^*:H_T^*(X^\gg) \ts_{H_T^*} \cO^h_{\cC_T}(U-\aa)
\ra H_T^*(X^\gg) \ts_{H_T^*} \cO^h_{\cC_T}(U-\bb)$.  But since
$\bb\prec\gg$ it follows that $H(\bb)$ fixes $X^\gg$,
so  we can replace $X^\bb$ by $X^\gg$ in Corollary~\ref{phi-translation},
and obtain a map which we also denote by $\tau_{\bb-\aa}^*$.

Since $H=H(\bb)$ fixes $X^\bb$ and $X^\gg$, we can use
Proposition~\ref{fixed} to replace (to pick an example)
$H_T^*(X^\bb) \ts_{H_T^*} \cO^h_{\cC_T}(U-\bb)$ by
$H_K^*(X^\bb) \ts_{H_K^*} \cO^h_{\cC_T}(U-\bb)$, where $K=T/H$.
We have to show that the following diagram is commutative:
$$\xymatrix{
   H_K^*(X^\bb) \ts_{H_K^*} \cO^h_{\cC_T}(U-\aa)
              \ar[d]^{1\ts t_{\bb-\aa}^*} \ar[dr]^{i^*\ts 1}  &    \\
   H_K^*(X^\bb) \ts_{H_K^*} \cO^h_{\cC_T}(U-\bb)
              \ar[d]^{i^*\ts 1}                   &
        H_K^*(X^\gg) \ts_{H_K^*} \cO^h_{\cC_T}(U-\aa)
              \ar[dl]^{1\ts t_{\bb-\aa}^*}                 \\
   H_K^*(X^\gg) \ts_{H_K^*} \cO^h_{\cC_T}(U-\bb)    &
            }$$
Now observe that both composites are equal to $i^*\ts t_{\bb-\aa}^*$.
\ep

\bdefn \label{kt-sheaf}
Let $X$ be a finite $T$-$CW$ complex, and $\cU=(U_\aa)_{\aa\in\cC_T}$
an adapted cover of $\cC_T$.  We then define a sheaf
$\cF=\cK_T^*(X)$ on $\cC_T$ by gluing the sheaves $\cF_\aa$
defined in~\ref{local-sheaf} via the gluing maps $\phi_{\aa\bb}$
defined in~\ref{gluing-map}.
\edefn

Notice that we can remove the dependence of $\cK_T^*(X)$ on the
cover $\cU$ as follows: Let $\cU$ and $\cV$ be two covers adapted
to $X$.  Then any common refinement $\cW$ is going to be adapted
as well, and the corresponding maps of sheaves
$\cF_\cU \ra \cF_\cW \la \cF_\cV$ are isomorphisms
on stalks, hence isomorphisms of sheaves.  Therefore we can omit the
subscript $\cU$.

\bthm \label{coh-thy}
$\cK_T^*(\blank)$ is a $T$-equivariant cohomology theory on finite
complexes, with values in the category of coherent holomorphic
sheaves of $\Z/2$-graded $\cO^h_{\cC_T}$-algebras.
\ethm

\bp
We start by showing that the construction of $\cK_T^*(X)$ is natural.
Let $X$ and $Y$ be two finite $T$-$CW$ complexes, and $f:X\ra Y$ a
$T$-equivariant map between them.  We want to define a map of
sheaves $f^*:\cK_T^*(Y)\ra \cK_T^*(X)$ with the properties
that $1_X^* = 1_{\cK_T^*(X)}$ and $(fg)^*=g^* f^*$.  Consider two
$T$-cell decompositions of $X$ and $Y$, and let $\cA$ be the collection
of all compact subgroups $L$ of $T$ such that $\D^n\times (T/L)$ appears
in the cell decomposition of either $X$ or $Y$.  Let $\cU$ be a cover
adapted to $\cA$.  Since $f$ is $T$-equivariant, for each $\aa$
we get by restriction a map $f_\aa : X^\aa\ra Y^\aa$.  This induces
for all $U\subseteq U_\aa$ a map
$f_\aa^* \ts 1: H_T^*(Y^\aa) \ts_{H_T^*} \cO^h_{\cC_T}(U-\aa)
 \ra H_T^*(X^\aa)\ts_{H_T^*}\cO^h_{\cC_T}(U-\aa)$, which commutes with
restrictions.  Therefore, we obtain a sheaf map
$f_\aa^*:\cF_\aa(Y) \ra \cF_\aa(X)$.  To get our global map $f^*$,
we only have to check that the maps $f^*_\aa$ glue well, i.e.\ that
they commute with the gluing maps $\phi_{\aa\bb}$.  This follows easily
from the naturality of ordinary equivariant cohomology, and from the
naturality in $X$ of the isomorphism 
$H_T^*(X^\bb) \cong H_K^*(X^\bb) \ts_{H_K^*} H_T^*$.  (See 
Proposition~\ref{fixed} with $Z=X^\bb$, $H=H(\bb)$ and $K=T/H$.)

Let $(X,A)$ be a pair of finite $T$-$CW$ complexes, i.e.\ $A$ is a closed
subspace of $X$ and the inclusion map $A\ra X$ is $T$-equivariant.
We then define $\cK_T^*(X,A)$ as the kernel of the map
$j^*: \cK_T^*(X/A) \ra \cK_T^*(point)$, where $j:point=A/A \ra X/A$ is
the inclusion map.  If $f:(X,A)\ra (Y,B)$ is a map of pairs of
finite $T$-$CW$ complexes, then $f^*:\cK_T^*(Y,B)\ra \cK_T^*(X,A)$ is
defined as the unique map induced on the corresponding kernels from
$f^*:\cK_T^*(Y)\ra \cK_T^*(X)$.

The last definition we need is of the coboundary map.  If $(X,A)$
is a pair of finite $T$-$CW$ complexes, we want to define
$\dd:\cK_T^*(A) \ra \cK_T^{*+1}(X,A)$.  This is obtained by
gluing the maps
 $$\dd_\aa \ts 1: H_T^*(A^\aa)\ts_{H_T^*}\cO^h_{\cC_T}(U-\aa) \ra
   H_T^{*+1}(X^\aa,A^\aa) \ts_{H_T^*} \cO^h_{\cC_T}(U-\aa) \comma$$
where $\dd_\aa:H_T^*(A^\aa) \ra H_T^{*+1}(X^\aa,A^\aa)$ is the usual
coboundary map.  The maps $\dd_\aa\ts 1$ glue well, because $\dd_\aa$
is natural.

To check the usual axioms of a cohomology theory (naturality, exact
sequence of a pair, and excision) for $\cK_T^*(\blank)$, recall
that it was obtained by gluing the sheaves $\cF_\aa$ along the maps
$\phi_{\aa\bb}$.  Since the sheaves $\cF_\aa$ were defined using
$H_T^*(X^\aa)$, the properties of ordinary $T$-equivariant
cohomology pass on to $\cK_T^*(\blank)$, as long as tensoring
with $\cO^h_{\cC_T}(U-\aa)$ over $H_T^*$ preserves exactness.
But this is implied by Proposition~\ref{flat-completion}.
\ep

Now, $\cC_T=\specm K_T^*$ is a nonsingular affine complex variety,
so it is Stein.  By a generalization of Theorem B of Cartan, the
sheaf cohomology vanishes on $\cC_T$ in positive dimensions for any
coherent sheaf.  (See Gunning and Rossi~\cite{GR}.)  This implies
that $\Gg$, the global sections functor, is exact, so we get the
following result:

\bcor \label{coh-thy-sections}
$\Gg\cK_T^*(\blank)$ is an $T$-equivariant cohomology theory 
on finite complexes with values in the category of $\Z/2$-graded
algebras.
\ecor

We want to investigate a few more useful properties of $\cK_T^*(\blank)$.

\bprop \label{point}
The value of the theory on a point is given by
$\cK_T^* \cong \cO^h_{\cC_T}$.
\eprop

\bp
Notice that if $X$ is a point, translation by $\aa$ induces an
isomorphism of the corresponding sheaf $\cF_\aa$ on $U_\aa$
to $\cO^h_{\cC_T \; \mid U_\aa}$.  Via this isomorphism, the gluing map
$\phi_{\aa\bb}$ is the usual restriction of holomorphic functions.
Therefore $\cK_T^*\cong\cO^h_{\cC_T}$.
\ep

\bprop \label{fflat}
Let $\Gg$ be the global sections functor.  Then the ring
$\Gg\cK_T^*$ is faithfully flat over $K_T^*$.
\eprop

\bp
Recall that we defined $\cC_T$ as the complex algebraic variety
$\specm K_T^*$.  As in the proof of Proposition~\ref{flat-completion},
denote by $\cO$ the algebraic structure sheaf of $\cC_T$, and by
$\cO^h$ the analytic structure sheaf.  Let $X=\cC_T$.  We need to
show that the natural map $\cO(X) \ra \cO^h(X)$ is faithfully
flat.  We first show it is flat:  By the local characterization of
flatness (e.g.\ Theorem (3.J) in Matsumura~\cite{Ma}), it is enough
to show that for any $x\in U$, the natural map $\cO_x \ra \cO^h_x$ is
flat.  But this is true, and the proof goes exactly as in
Proposition~\ref{flat-completion}.

To prove that $\cO(X) \ra \cO^h(X)$ is in fact faitfully flat, we
also need to show that the induced map of spectra of maximal ideals, 
$\specm \cO^h(X) \ra \specm \cO(X)$ is surjective (see 
Theorem 3 (4.D) in~\cite{Ma}).  This is clearly true, since maximal
ideals of $\cO(X)$ and $\cO^h(X)$ correspond to the points $x\in X$.
\ep

Since $\cK_T^*(X)$ was obtained from gluing the sheaves $\cF_\aa$,
its stalks are the same as those of $\cF_\aa$.

\bprop \label{stalk}
If $X$ is a finite $T$-$CW$ complex and $\aa\in\cC_T$ , then the stalk
of the $\cK_T^*(X)$ at $\aa$ is
    $$\cK_T^*(X)_\aa = H_T^*(X^\aa) \ts_{H_T^*}\cO^h_{\cC_T,0}\period$$
\eprop

\bprop \label{product}
Let $A$ and $B$ are compact abelian Lie groups, $X$ is a finite
$A$-$CW$ complex, and $Y$ is a finite $B$-$CW$ complex.  Then
 $$\Gg\cK_{A\times B}^*(X\times Y) \cong
     \Gg\cK_A^*(X) \ts_\C \Gg\cK_B^*(Y) \period$$
\eprop

\bp \label{holo-K-theory}
This follows from $\cC_{A\times B}\cong \cC_A \times\cC_B$ and from
$H_{A\times B}^*(X\times Y) \cong H_A^*(X) \ts_\C H_B^*(Y)$, where
the last isomorphism comes from the K\"{u}nneth formula.
\ep

We now return to the usual equivariant $K$-theory, and define
its holomorphic extension
    $$ \ko_T^*(X) = K_T^*(X) \ts_{K_T^*} \Gg\cK_T^* \period$$
Since $K_T^*(\blank)$ is a $T$-equivariant cohomology theory,
it follows from Lemma~\ref{fflat} that $\ko_T^*(\blank)$ is
a $T$-equivariant cohomology theory, with values in $\Z/2$-graded
algebras.

\section{The equivariant Chern character}

Let $T$ be a compact abelian Lie group and $X$ a finite $T$-$CW$
complex.  In this section we construct a natural isomorphism
    $$CH_T: \ko_T^*(X) \ra \Gg\cK_T^*(X) \period$$
The main ingredient in the definition of $CH_T$ is the construction
of a natural ring map $CH_T: K_T^0(X,\Z) \ra \Gg\cK_T^*(X)$.  
Let $E\ra X$ be a complex $T$-vector bundle.  We want $CH_T(E)$ to
be a section in $\cK_T^*(X)$, so we need to know its germ at a point
$\aa\in\cC_T$.  The germ $CH_T(E)_\aa$ should be an element of
$\ho_T^*(X^\aa)$, so we would be quite tempted to define it as
the equivariant Chern character $ch_T$ of the restriction of $E$ to
$X^\aa$ (we will define the equivariant Chern character in the next
subsection).  However, there are a few problems.  The first one is
that in order to get $ch_T$ to take values in $\ho_T^*(\blank)$,
we need  some sort of holomorphicity result on $ch_T$.  This is
done in Lemma~\ref{cht-holo}.  The second problem is that when
we vary $\aa\in\cC_T$, the equivariant Chern character does not glue
well.  To solve this problem, we modify the definition of $CH_T$
by twisting $ch_T$ accordingly: see Definition~\ref{cht-twisted}.

\vs
\bc
{\sc 3.1. Preliminaries}
\ec
\vs

Let $E\ra X$ be a complex $T$-vector bundle.  Then the equivariant
Chern character $ch_T(E)$ can be defined as $ch(E_T)$, where $ch$
is the usual Chern character, and $E_T \ra X_T$ is the Borel
construction.  This a priori is an element of $H_T^{**}(X)$,
because the Chern character is defined using the power series
$e^x$ (see the Appendix of Rosu~\cite{Ro} for a discussion of
equivariant characteristic classes).  So we get a map
   $$ ch_T: K_T^*(X) \ra H_T^{**}(X) \period$$
Also, denote the $n$-th Chern class of $E_T$ by $c_n(E)_T$.
This belongs to the noncompleted ring $H_T^*(X)$.

\blemma \label{cht-holo}
Let $E\ra X$ be a complex $T$-vector bundle over a finite
$T$-$CW$ complex.   Then its $T$-equivariant Chern character,
$ch_T(E)$, which a priori is an element of $H_T^{**}(X)$,
actually lies in $\ho_T^*(X)$.
\elemma

\bp
First assume that $E$ is a line bundle.  Then from the definition of
the Chern character we get $ch_T(E)=e^{c_1(E)_T}$.  We have already
seen that $c_1(E)_T \in H_T^*(X)$.  Since $X$ is finite,
$H_T^*(X)$ is a finite module over $H_T^*$.  Let $a_1,\ldots,a_m$ be
a set of generators for $H_T^*(X)$.  Choose a set of elements
$f^k_{ij}\in H_T^*$ so that
     $$a_i\cdot a_j = \sum_{k=1}^m f^k_{ij} a_k \hspace{3mm}
              \forall\; i,j,k \in 1,\ldots,m \period$$
Denote by $c=c_1(E)_T$.  The element $c$ is in $H_T^*(X)$, so there
exist elements $g_i \in H_T^*$ such that
    $$c= g_1 a_1 + \cdots + g_m a_m \period$$
We can also calculate $c^2= \sum_k (\sum_{i,j} g_i g_j f^k_{ij}) a_k$,
$c^3=\sum_p (\sum_{i,j,k,l}g_i g_j g_l f^k_{ij} f^p_{kl}) a_p$, etc.

Choose some coordinates on $t_\C$ such that $H_T^*=\C[u_1,\ldots,u_p]$.
We say that a polynomial $\phi\in \C[u_1,\ldots,u_p]$ is dominated by
another polynomial $\psi$ if all coefficients of $\psi$ are positive
and all the coefficients of $\phi$ are less in absolute value than the
coefficients of $\psi$.  Also, if $\mu, \nu \in H_T^*(X)$ then
we can write $\mu = \mu_1 a_1 +\cdots + \mu_m a_m$ and
$\nu = \nu_1 a_1 + \cdots + \nu_m a_m$, with $\mu_i,\nu_i \in H_T^*$.
We then say that $\mu$ is dominated by $\nu$ if $\mu_i$ is dominated by
$\nu_i$ for all $i$.  Similar definitions of domination 
can be made for the power series ring $H_T^{**}$ and for $H_T^{**}(X)$.

Let $\ll \in H_T^*=\C[u_1,\ldots,u_p]$ be a polynomial such that
$f^k_{ij}$ and $g_i$ are dominated by $\ll$ for all $i,j,k$.  Then
one can show by induction that $c^n$ is dominated by 
$(2n-2)\ll^{2n-1}(a_1+\cdots+a_m)$, which in turn is dominated by
$2n\ll^{2n-1}(a_1+\cdots+a_m)$.  So $e^c$ is dominated by
$1+\sum_{n\geq 1} \tfrac2{(n-1)!}\, \ll^{2n-1} (a_1+\cdots+a_m) =
1+\sum_{n\geq 0} \tfrac2{n!} \,\ll^{2n+1} (a_1+\cdots+a_m) = 
1+2\ll e^{2\ll} (a_1+\cdots+a_m)$.  But this last element belongs to
$H_T^*(X) \ts_{H_T^*} \cO^h_{t_\C,0} = \ho_T^*(X)$, so $e^c$ also belongs
to $\ho_T^*(X)$.

Second, assume $E$ is a rank $n$ complex $T$-vector bundle.
Let $S_E$ be the splitting space of $E$ (or the flag variety of
$E$---see Bott and Tu~\cite{BTu}).  By construction $S_E$ comes with an
equivariant map $\pi:S_E \ra X$.  Then the splitting principle says
that the pull-back bundle $\pi^*(E)$ decomposes as a sum of
$T$-equivariant line bundles over $S_E$.  Say 
$\pi^*(E)\cong L_1 \oplus \cdots \oplus L_n$.  Calculation using the
Leray--Serre spectral sequence shows that
  $$H_T^*(S_E)= 
   H_T^*(X)[x_1,\ldots,x_n]\;/\langle \ss_i(x_1,\ldots,x_n)-
                            c_i(E)_T\rangle \comma$$
where $\ss_i(x_1,\ldots,x_n)$ is the $i$'th symmetric polynomial in the
$x_j$'s.  The classes $x_j=c_1(L_j)_T$ are called the Chern roots of $E$.
Moreover, we can identify $H_T^*(X)$ as the subring of $H_T^*(S_E)$
generated by the polynomials in $H_T^*(X)[x_1,\ldots,x_n]$ which are
symmetric in the $x_j$'s.  By tensoring with $\cO^h_{t_\C,0}$ or $H_T^{**}$
the same statement is true about $\ho_T^*(X)$ and $H_T^{**}(X)$.

Now consider $ch_T(E) = e^{x_1} \cdots e^{x_n}$.  Since
$L_j$ is a line bundle and $x_j=c_1(L_j)_T$, the first part of the proof
implies that $e^{x_j}\in \ho_T^*(S_E)$ for all $j$.  Therefore
$ch_T(E) \in \ho_T^*(S_E)$, and since it is symmetric in the $x_j$'s
it follows that $ch_T(E) \in \ho_T^*(X)$, which is what we wanted.
\ep

We have just proved that $ch_T(E)$ is the germ of a holomorphic
class, i.e.\ an element of $\ho_T^*(X)$.  By looking more carefully
at the preceding proof, one can see in fact that we proved 
a stronger result:

\bcor \label{cht-global-holo}
With the same notations as in Lemma~\ref{cht-holo}, $ch_T(E)$ is
a global holomorphic class, i.e.\ an element of
$H_T^*(X) \ts_{H_T^*} \Gg\cO^h_{\cC_T}$.
\ecor

Now, if we extend $ch_T(E)$ on a small neighborhood $U$ of
$0\in\cC_T$, we can regard it as an element of
$H_T^*(X) \ts_{H_T^*} \cO^h_{\cC_T}(U)$.  It is important to see
what happens to $ch_T(E)$ when it is translated by the map
$\tau_{\bb-\aa}^*$ from  Proposition~\ref{phi-translation}.

The basic case is when  $X$ is a point and $E$ is given by a
representation $V_\ll$ of $T$, with $\ll\in\T$.
Recall that $\cC_T=\Hom(\T,\cC)$, and consider $\aa\in\cC_T$.
Then we translate $ch_T(V_\ll)$ via the map
$\tau_\aa^*=t_\aa^*:\cO^h_{\cC_T}(U) \ra \cO^h_{\cC_T}(U+\aa)$.

\blemma \label{cht-translation}
Let $T$ be a compact Abelian group, and $T^0$ the connected component
containing the identity.  Let $\aa\in\cC_{T^0}$ and $\ll\in\T$.
Then, with the notations above,
   $$ t_\aa^*ch_T(V_\ll) = \aa(\ll)ch_T(V_\ll) \period$$
\elemma

\bp
The proof of this lemma is mostly formal, and just makes intensive
use of the identifications and definitions we have made so far.
Start by using the exponential map
  $$\exp: t_\C = \Ll_T\ts_\Z\C \ra 
        \Ll_T\ts_\Z\cC = \cC_{T^0} \ra \cC_T \period$$
The element $\aa$ is in the image, so pick $a\in t_\C$ such that
$\exp(a)=\aa$.  Denote by $l\in\Ll_T^*$ the element corresponding
to $\ll$ via the map $\T \ra \Ll_T^*$.  Then one can 
apply $l$ to $a\in\Ll_T\ts_\Z\C$ and get a complex number
that we denote by $l(a)$.  Now it is easy to check the formula
$\aa(\ll)=\exp(l(a))$.  We also know that 
$ch_T(V_\ll)=\exp\bigl(c_1(V_\ll)\bigr)$.  So, via the exponential
map, what we have to prove becomes
  $$ t_a^* c_1(V_\ll) = l(a) + c_1(V_\ll) \comma$$
with the equality being now regarded in $H_T^*$.
Let us look more closely at $c_1(V_\ll)$.  We saw in the proof
of Proposition~\ref{polynomial} that there is an identification
$H_T^*=S(t_\C^*)$, and that the class $c_1(V_\ll)\in H_T^*$
can be identified to $l$ if this is regarded in $S(t_\C^*)$ via
$\Ll_T^* \subset t_\C^* \subset S(t_\C^*)$.  Denote by $l(\blank)$
the polynomial function in $S(t_\C^*)$ corresponding to $l$.
Then we have to prove that
     $$ t_a^* l(\blank) = l(a) + l(\blank) \period$$
But this is obvious, it is just saying that $l(\blank)$ is a
linear function.
\ep

\vs
\bc
{\sc 3.2. Construction of} $CH_T$ 
\ec
\vs

We define a multiplicative natural map
  $$CH_T: K_T^0(X,\Z) \ra \Gg\cK_T^*(X) \period$$
Let $E\ra X$ be a complex $T$-vector bundle.  Let $\aa\in\cC_T$, and
denote by $H=H(\aa)$.  Then $\aa \in \cC_H$.
By Proposition~\ref{ct}, $\cC_H \cong \Hom_\Z(\H,\cC)$, so we can
think of $\aa$ as a group map $\aa: \H \ra \cC$.  The space
$X^\aa$ has a trivial action of $H$, so the restriction $E_{\mid X^\aa}$
of $E$ to $X^\aa$ has a fiberwise decomposition by irreducible
characters of $H$:
   $$E_{\mid X^\aa} \cong \oplus_{\ll\in \H} E(\ll) \comma$$
where $E(\ll)$ is the $T$-vector bundle where $h\in H$ acts by complex
multiplication with $\ll(h)$.

It would be tempting to define the germ of $CH_T(E)$ at $\aa$ to be
$ch_T(E_{\mid X^\aa})$, but these germs would not glue well to give a
global section of $\cK_T^*(X)$.  Instead, we do the following:

\bdefn \label{cht-twisted}
Let $\aa\in\cC_T$ and $H=H(\aa)$.  Then the germ of
$CH_T(E)$ at $\aa$ is defined to be
   $$CH_T(E)_\aa = \sum_{\ll\in\H} \aa(\ll) ch_T E(\ll) \period$$
\edefn

\bprop \label{cht-glue}
The germs $CH_T(E)_\aa$ glue to a global section
$CH_T(E)\in \Gg\cK_T^*(X)$.
\eprop

\bp
We notice that, by Lemma~\ref{cht-holo}, $CH_T(E)_\aa$ does indeed
belong to $\ho_T^*(X^\aa)$, which by Proposition~\ref{stalk} is the
stalk of $\cK_T^*(X)$ at $\aa$.  Fix $(U_\aa)_{\aa\in\cC_T}$
a cover of $\cC_T$ adapted to $X$.

Let $\aa,\bb\in\cC_T$ with $U_\aa\cap U_\bb\neq\emptyset$ and
$\aa\prec\bb$.  This implies $\aa,\bb \in \cC_{H(\bb)}$ and also
$H(\aa) \subseteq H(\bb)$.  Denote by $L=H(\aa)$ and $H=H(\bb)$.
Condition 4 of Definition~\ref{adapted-defn} implies that
$\bb-\aa \in \cC_{H^0}$.
Now we have to prove that $\phi_{\aa\bb}CH_T(E)_\aa =CH_T(E)_\bb$,
i.e.\ that 
$\tau_{\bb-\aa}^* \sum_{\ll\in\L} \aa(\ll) ch_T E(\ll)_{\mid X^\bb} = 
          \sum_{\mu\in\H} \bb(\mu) ch_T E(\mu)$.
Consider the surjective map $j: \H \ra \L$ induced by
the inclusion $L \inc H$.  If $\ll\in \L$, we have
$E(\ll)_{\mid X^\bb} = \sum_{\mu\in j^{-1}(\ll)} E(\ll)$.  Therefore
it is enough to show that for all $\mu\in \H$ we have
$\tau_{\bb-\aa}^* \bigl(\aa(\ll) ch_T E(\mu)\bigr) = \bb(\mu)ch_T E(\mu)$,
where $\ll=j(\mu)$.  But this is equivalent to
$\tau_{\bb-\aa}^* ch_T E(\mu) = (\bb-\aa)(\mu) ch_T E(\mu)$.
Denote by $\gg=\bb-\aa \in \cC_{H^0}$.  So it is enough to show
that, for all $\gg \in \cC_{H^0}$ and $\mu\in\H$,
 $$\tau_\gg^* ch_T E(\mu) = \gg(\mu) ch_T E(\mu) \period$$

Let $K=T/H$ and $Y=X^\bb$.  Proposition~\ref{fixed} applied to
equivariant $K$-theory gives a natural isomorphism
  $$ K_K^*(Y) \ts_{K_K^*} K_T^* \iso  K_T^*(Y) \period$$
Via the identification above, we can think of $E(\mu)$ as a tensor
product $F\ts V(\mu')$, with $F$ a $K$-bundle and $\mu'\in\T$ some
element in the preimage of $\mu$ via the map $\T\ra\H$.
(At least, we know that $E(\mu)$ is generated by such elements.)
Via the same identification, translation by $\gg\in\cC_{H^0}$
becomes
   $$\tau_\gg^* \mapsto \tau_{\pi(\gg)}^* \ts \tau_\gg^* \comma$$
where $\pi(\tau)$ is the image of $\gg$ via the natural map
$\pi:\cC_T \ra \cC_K$, and the second $\gg$ is regarded in $\cC_T$
via the usual inclusion $\cC_{H^0}\ra\cC_T$.  But notice that
$\pi(\tau)=0$, because of the exact sequence
$0 \ra \cC_H \ra \cC_T \ra \cC_K \ra 0$.  Also, $ch_T E(\mu)$ becomes
$ch_K F \ts ch_T V(\mu') \in \ho_K^*(Y) \ts_{\ho_K^*} \ho_T^*$.
So via the above correspondence we have
  $$ \tau_\gg^* ch_T E(\mu) \mapsto
            ch_K F \ts \tau_\gg^* ch_T V(\mu') \period$$
Since $\gg\in\cC_{H^0}$, it follows that $\gg\in\cC_{T^0}$, and
it is sufficient for us to show that, for all $\gg \in \cC_{T^0}$
and $\mu'\in\T$,
  $$\tau_\gg^* ch_T V(\mu') = \gg(\mu) ch_T V(\mu') \period$$
But this is directly implied by Lemma~\ref{cht-translation},
so we are done.
\ep

We have just finished constructing a natural map
$CH_T: K_T^0(X,\Z) \ra \Gg\cK_T^*(X)$.  By taking the suspension of
$X$ instead of $X$, this induces a map
$CH_T: K_T^*(X,\Z) \ra \Gg\cK_T^*(X)$. 
One can check easily that $CH_T$ is a ring map, since $ch_T$ is.
Because $\Gg\cK_T^*(X)$ is a $\C$-algebra and $CH_T$ is a ring map,
we can now extend $CH_T$ to a natural map of $\C$-algebras
$CH_T: K_T^*(X) \ra \Gg\cK_T^*(X)$.  Finally, making $X$ a point,
we get a ring map $K_T^* \inc \Gg\cK_T^*$, so if we
extend $CH_T$ by this, we obtain the desired natural map
    $$CH_T: \ko_T^*(X) \ra \Gg\cK_T^*(X) \period$$

\bthm
$CH_T$ is an isomorphism of $T$-equivariant cohomology theories.
\ethm

\bp
Because of the Mayer--Vietoris sequence, it is enough to verify the
isomorphism for ``equivariant points'' of the form $T/L$, with $L$ a
compact subgroup of $T$.  Choose an identification
 $$T= (S^1)^p\times (S^1)^q \times \prod_{i=1}^r \Z_{m_i}$$
such that via this identification
 $$L= (S^1)^p\times \prod_{j=1}^q \Z_{n_i} 
              \times \prod_{i=1}^r \Z_{l_i} \period$$
Then $T/L= (1)^p \times \prod_{j=1}^q (S^1/\Z_{n_i})
              \times \prod_{i=1}^r (\Z_{m_i}/\Z_{l_i})$.

We use now Proposition~\ref{product}, which is also true if we
replace $\Gg\cK$ by $\ko$-theory (because it is true for $K$-theory).
Since the map $CH$ commutes
with the isomorphisms of Proposition~\ref{product}, it is enough
to check that the following maps are isomorphisms:
\bi
\item[(a)] $CH_{S^1}:\ko_{S^1}^* \ra \Gg\cK_{S^1}^*$;
\item[(b)] $CH_{S^1}:\ko_{S^1}^*(S^1/\Z_n) \ra \Gg\cK_{S^1}^*(S^1/\Z_n)$;
\item[(c)] $CH_{\Z_m}: \ko_{\Z_m}^*(\Z_m/\Z_l) \ra
       \Gg\cK_{\Z_m}^*(\Z_m/\Z_l)$.
\ei

To prove $(a)$, notice that $\cC_{S^1}=\cC$.  Then we have
$\ko_{S^1}^* = K_{S^1}^* \ts_{K_{S^1}^*} \Gg \cO^h_{\cC} = \Gg \cO^h_{\cC}$.
By Proposition~\ref{point}, $\Gg\cK_{S^1}^*=\Gg\cO^h_{\cC}$.  Now notice
that, by definition, the map $CH_{S^1}$ is the identity.

\vs

For $(b)$, denote $X=S^1/\Z_n$.  Then we have
$\ko_{S^1}^*(X)=K_{S^1}^*(X) \ts_{K_{S^1}^*} \Gg\cK_{S^1}^*=
K_{\Z_n}^* \ts_{K_{S^1}^*} \Gg\cO^h_{\cC}$.  But we know that
$K_{S^1}^*=\C[z^{\pm 1}]$ and
$K_{\Z_n}^*=\C[z^{\pm 1}]/ \langle z^n-1\rangle$.
So we deduce $\ko_{S^1}^*(X) = \Gg\cO^h_{\cC}/ \langle z^n-1\rangle$.
This last ring can be identified with
$\C[z^{\pm 1}]/ \langle z^n-1\rangle$, since the condition
$z^n=1$ makes all power series finite.  In conclusion,
$\ko_{S^1}^*(X) = K_{\Z_n}^*= \C[z^{\pm 1}]/ \langle z^n-1\rangle$.

Let us now describe the sheaf $\cF=\cK_{S^1}^*(X)$.  Let $\aa\in\cC$.
If $\aa\notin\Z_n$, $X^\aa=\emptyset$, so the stalk of $\cF$ at $\aa$
is zero.  If $\aa\in\Z_n$, $X^\aa=X$, and the stalk of $\cF$ at $\aa$
is $H_{\Z_n}^* \ts_{H_{S^1}^*} \cO^h_{\C,0}$.  But $H_{\Z_n}^* = \C$,
concentrated in degree zero ($H_{\Z_n}^*$ is $\Z$-torsion in higher
degrees, so the components in higher degree disappear when we tensor
with $\C$).  It follows that
$\cF$ is a sheaf concentrated at the elements of $\Z_n$, where it has
the stalk equal to $\C$.  Then the global sections of $\cF$ are
$\Gg\cK_{S^1}^*(X)=\C\oplus \cdots \oplus \C$, $n$ copies, one for
each element of $\Z_n$.

The map $CH_{S^1}:\ko_{S^1}^*(X) \ra \Gg\cK_{S^1}^*(X)$ comes from
the ring map
$CH_T: \C[z^{\pm 1}]/ \langle z^n-1\rangle \ra \C\oplus \cdots \oplus \C$.
Since $z$ generates the domain of $CH_T$, it is enough to see where $z$
is sent.  Let $\Z_n=\{1,\eps,\eps^2,\ldots,\eps^{n-1}\}$.  Then
$z$ represents the standard irreducible representation $V=V(\eps)=\C$
of $\Z_n$, where $\eps$ acts on $\C$ by complex multiplication with
$\eps$, which is regarded as an element of $\cC$.  Notice that
$V$ corresponds to the element $\ll=\eps\in \widehat{\Z_n}=\Z_n$.
$c_1(V)_{S^1}=0$, because $c_1(V)_{S^1}$ lies in $H_{\Z_n}^2=0$.
Then $ch_{S^1}(V)=e^{c_1(V)_{S^1}}=e^0=1$, and
the stalk of $CH_{S^1}(V)$ at $\aa\in\Z_n$ is
$CH_{S^1}(V)_\aa = \aa(\eps) = \aa$.  Therefore $CH_{S^1}$ sends
$z$ to $(1,\eps,\eps^2,\ldots,\eps^{n-1}) \in \C\oplus \cdots \oplus \C$.
One can easily check that this map is an isomorphism.

\vskip2mm

For $(c)$, denote $X=\Z_m/\Z_l$.  As in $(b)$, 
$\ko_{\Z_m}^*(X) =  K_{\Z_l}^*=\C[z^{\pm 1}]/ \langle z^l-1\rangle$, and
$\Gg\cK_{\Z_m}^*(X)=\C\oplus \cdots \oplus \C$, $l$ copies.
The proof that $CH_{\Z_m}$ is an isomorphism is the same
as above.
\ep

\appendix
 
\section{Applications}

We now give applications of the construction in this paper.  First
we use the Chang--Skjelbred theorem in equivariant cohomology to infer
the corresponding result for equivariant $K$-theory.  Then as a corollary
we show how a result about the equivariant cohomology of GKM manifolds
can be extended to equivariant $K$-theory.  Along the way we
need a natural splitting that does not seem to have been 
noticed before in this area.

\bdefn
Let $X$ be a compact $T$-manifold, for $T$ a compact abelian Lie group.
We say that $X$ is equivariantly formal if the equvariant
cohomology spectral sequence collapses at the $E_2$ term.
\edefn

Many interesting $T$-spaces are equivariantly formal; for example any
subvariety of complex projective space preserved by a linear action,
or symplectic manifold with a Hamiltonian action.
Our reference for equivariantly formal spaces is~\cite{GKM}.
We need three results about them:  The first is that the map
$H^*_T(X) \to H^*_T(X^T)$ is an injection.  The second is that for any
$H$ compact subgroup of $T$, $X^H$ is also equivariantly formal.  The
third is due to Chang and Skjelbred~\cite{CS} (see also~\cite{BV} for
a proof):

\bthm \label{chskj}
Let $X$ be an equivariantly formal $T$-manifold, and let $X_1$ be its
equivariant 1-skeleton, i.e.\ $X_1$ is the set of points in $X$ with
stabilizer of codimension at most one.  We have inclusion maps
$i:X^T\ra X$ and $j:X^T\ra X_1$.  
Then the map $i^*:H_T^*(X) \ra H_T^*(X^T)$ is injective, and 
the maps $i^*$ and $j^*:H_T^*(X_1) \ra H_T^*(X^T)$ have the same image.  
\ethm

The ring $H_T^*(X_1)$, in the notation of the above theorem, has
not received much study. It is typically much bigger
than $H_T^*(X)$, and though $H_T^*(X)$ injects into it, it does not
inject into $H_T^*(X^T)$. These phenomena can be seen in the case
of $T^2$ acting on $X=\C\sP^2$, where $X^1$ is a cycle of three $\C\sP^1$'s
and therefore has $H^1$, not seen in either $H_T^*(X)$ or $H_T^*(X^T)$.

\blemma\label{splitting}
In the notation of Theorem \ref{chskj},
there is a natural identification $H_T^*(X_1) = H_T^*(X) \oplus \ker j^*$.
\elemma

\bp
By Theorem~\ref{chskj} the images of the two maps
$i^*:H_T^*(X) \ra H_T^*(X^T)$ and $j^*:H_T^*(X_1) \ra H_T^*(X^T)$ are
the same.  But $i^*$ is injective, so we can identify $H_T^*(X)$ with
the image of $i^*$.  This implies that $j^*$ factors through a map
$H_T^*(X_1) \ra H_T^*(X)$, and this yields a splitting
$H_T^*(X_1) = H_T^*(X) \oplus \ker j^*$.
\ep

This natural splitting sheafifies, allowing us to extend both results
to $K$-theory. 

\bthm \label{gkm-ktheory}
We use the same notations as in Theorem~\ref{chskj}.  Then
$i^*:K_T^*(X) \ra K_T^*(X^T)$ is injective, and the maps
$i^*$ and $j^*:K_T^*(X_1) \ra K_T^*(X^T)$ have the same image.
\ethm

\bp
Let $\aa\in\cC_T$.  Any compact $T$-manifold admits a decomposition as
a finite $T$-$CW$ complex (see for example Allday and Puppe~\cite{AP}).
Let $X=\bigcup_i \D^{n_i}\times (T/H_i)$ be such a cell decomposition.
We saw that if $K$ is a compact subgroup of $T$,
$X^K=\sum_{i:K\subseteq H_i} \D^{n_i}\times (T/H_i)$.  In particular,
this implies that $(X_1)^\aa = (X^\aa)_1$.

Let $Y=X^\aa$, which is again equivariantly formal.  By Lemma \ref{splitting}
there is a natural identification $H_T^*(Y_1) = H_T^*(Y) \oplus \ker j^*$.
By Proposition~\ref{flat-completion} the map $H_T^*\ra \ho_T^*$ is
flat, so tensoring with $\ho_T^*$ over $H_T^*$ yields a splitting
$\ho_T^*(Y_1) = \ho_T^*(Y) \oplus \ker j^*$.  Now, we observed above
that $(Y_1)^\aa = (Y^\aa)_1$.  So we finally get a splitting
$\ho_T^* \bigl((X_1)^\aa\bigr) = \ho_T^*(X^\aa) \oplus \ker j^*$.
This is compatible with the gluing maps of the sheaf $\cK_T^*(X)$, so
we get $\cK_T^*(X_1) = \cK_T^*(X) \oplus \ker j^*$.

The upshot of the above discussion is that
$i^*: \cK_T^*(X) \ra \cK_T^*(X^T)$ is injective (since it is injective
on stalks), and $\cK_T^*(X_1) = \cK_T^*(X) \oplus \ker j^*$.
The global section functor $\Gg$ is left exact,
so $i^*: \Gg \cK_T^*(X) \ra \Gg \cK_T^*(X^T)$ is injective and
$\Gg \cK_T^*(X_1) = \Gg \cK_T^*(X) \oplus \ker j^*$.  This implies
that $i^*$ and $j^*$ have the same image in $\Gg \cK_T^*(X^T)$,
namely $\Gg \cK_T^*(X)$.  (Notice we couldn't have done this without
using the splitting, because $\Gg$ is not right exact, so it doesn't
commute with the image functor.)

Now recall that we have a natural isomorphism
$CH_T:\ko_T^*(X) \ra \Gg\cK_T^*(X)$.  Translating
the above results via $CH_T$, we obtain that the maps
$i^* : \ko_T^*(X) \ra \ko_T^*(X^T)$ and $j^*:\ko_T^*(X_1) \ra \ko_T^*(X^T)$
have the same image.  But $\ko_T^*(X) = K_T^*(X) \ts_{K_T^*} \ko_T^*$
and by Lemma~\ref{fflat} the map $K_T^* \ra \ko_T^*$ is faithfully flat.
So we obtain that $i^*$ and $j^*$ have the same images in $K_T^*(X^T)$,
which is what we wanted.

\ep

An alternative proof of Theorem~\ref{gkm-ktheory} can be given by
noticing that the sheaf maps $i^*:\cK_T^*(X) \ra \cK_T^*(X^T)$ and
$j^*:\cK_T^*(X_1) \ra \cK_T^*(X^T)$ have the same image (because
they have the same image at the level of stalks).  But the global
section functor is exact, since we work with coherent sheaves over
Stein manifolds (see the comment before
Corollary~\ref{coh-thy-sections}).  It follows that the maps
$i^* : \ko_T^*(X) \ra \ko_T^*(X^T)$ and $j^*:\ko_T^*(X_1) \ra \ko_T^*(X^T)$
have the same image, and the proof proceeds as before.

In \cite{GKM} a special case of this is studied, in which $X^T$ is
discrete and $X_1$ is a union of $S^2$'s; these are called 
\emph{balloon manifolds} or \emph{GKM manifolds}.
(An interesting example of GKM manifolds are toric varieties.)
In this case it is easy to calculate the image of restriction from
$H^*_T(X_1)$, by reducing it to the case of $H^*_T(S^2)$. The theorem
above lets us extend this result to $K$-theory.

\bcor
Let $X$ be a GKM manifold, and $i^* : K_T^*(X) \to K_T^*(X^T)$ the 
restriction map.  Then $i^*$ is an injection; and a class
$\aa\in K_T^*(X^T)$ is in the image if for each 2-sphere
$B \subseteq X_1$ with fixed points $N$ and $S$, the difference
$\aa|_N - \aa|_S \in K_T$ is a multiple of the K-theoretic Euler
class of the tangent space $T_N B$.
(Technically, we can only take the Euler class once we orient $T_N B$,
but either orientation leads to the same condition on $\aa$.)
\ecor

If $T = (S^1)^n$, we can identify $K_T^*$ with Laurent polynomials,
and this condition says that the difference $\aa|_S-\aa|_N$
of Laurent polynomials must be a multiple of $1-w$, where $w$ is the
weight of the action of $T$ on $T_N B$. (Again, we can only speak of
the weight $w$ once we orient this $\R^2$-bundle over the point $N$,
but it doesn't matter because being a multiple of $1-w$ is the same
as being a multiple of $1-w^{-1}$. In most examples one has around
a $T$-invariant almost complex structure with which to orient 
all these tangent spaces simultaneously.)

After finishing this paper, similar results in the algebraic case
have appeared in~\cite{VV}, particularly with regard to the extension
of Chang--Skjelbred's results to K-theory.  More specifically, their
Corollary 5.10 is the algebraic analogue of our Theorem~\ref{gkm-ktheory}.
Notice also that their results, remarkably, hold over $\Z$, while ours
only hold over $\C$.  We thank Angelo Vistoli for explaining his results
to us.

We also found out that Atiyah~\cite{At} proved a Chang-Skjelbred
lemma in the context of equivariant $K$-theory.  He not only did
it more or less at the same time as Chang and Skjelbred,
but his results are stronger:  Let $G$ be a torus and $X$ a compact
$G$-manifold.  Denote by $X_i$ the equivariant $i$-skeleton of $X$.
Then Atiyah shows that the long exact sequence of the pair
$(X\setminus X_i,X\setminus X_{i+1})$ splits into short exact sequences
$0 \lra K_G^{*-1}(X\setminus X_{i+1}) \llra{\dd} 
    K_G^*(X_{i+1}\setminus X_i) \lra K_G^*(X\setminus X_i) \lra 0$.
This in turn is equivalent to the having a long exact sequence
  $$ 0 \lra K_G^*(X) \lra K_G^*(X^0) \lra K_G^{*+1}(X_1\setminus X_0)
          \lra K_G^{*+2}(X_2\setminus X_1) \lra \cdots $$
As noted by Bredon~\cite{Br}, Atiyah's argument carries over to
equivariant cohomology with compact support and rational
coefficients.  In this context, the exactness up to the
term $K_G^{*+1}(X_1\setminus X_0)$ is just the assertion of the 
Chang--Skjelbred lemma.  In cohomological setting, the above
sequence may be thought of as the $E_1$ term of the spectral sequence
coming from the filtration of the Borel construction $X_G$ by the
subspaces $(X_i)_G$.  Moreover, its exactness is in fact equivalent
to the equivariant formality of $X$, i.e.\ to the freeness of 
$H_G^*(X)$.  (The direction Atiyah proved is the harder one.)
We thank Matthias Franz for pointing out these results to us.

\subsection{Acknowledgements}

We thank Victor Guillemin for persuading us to write this paper,
and Haynes Miller for constant guidance and support throughout
this project.  We also thank Michele Vergne for going carefully
through the paper and suggesting several corrections and improvements.
Finally, we thank Lars Hesselholt, Payman Kassaei and Behrang Noohi
for helpful discussions.

\vskip5mm

\noindent
\textsc{Massachusetts Institute of Technology, Cambridge, MA}\\
\textit{E-mail address:} \texttt{ioanid@math.mit.edu}\\
\textsc{University of California at Berkeley, CA}\\
\textit{E-mail address:} \texttt{allenk@math.berkeley.edu}

\end{document}